\documentclass[12pt]{amsart}

\usepackage{amsmath,amsthm,amsfonts,latexsym,amssymb,enumerate,color}
\usepackage{graphicx}
\usepackage{enumitem}

\newcommand\CC{{\mathbb C}}

\newcommand\DD{{\mathbb D}}

\newcommand\RR{{\mathbb R}}

\newcommand\TT{{\mathbb T}}

\def\beq{\begin{equation}}
\def\eeq{\end{equation}}
\newtheorem{thm}{Theorem}[section]
\newtheorem{prop}[thm]{Proposition}

\newtheorem{cor}[thm]{Corollary}
\newtheorem{rem}[thm]{Remark}

\newcommand\beginpf{\noindent {\bf Proof:} \quad}

\newcommand\re{\mathop{\rm Re}\nolimits}

\def\beginpf{\begin{proof}}
\def\endpf{\end{proof}}

\newcommand{\ds}{\displaystyle}

\renewcommand\phi{\varphi}

\DeclareMathOperator{\dom}{dom}

\newcommand\Hol{\mathop{\rm Hol}}
\newcommand\LL{\mathcal L}

\begin{document}
\title[Composition operators on  the halfplane]{Composition operators on function spaces on the halfplane: 
spectra and semigroups}

\author{I. Chalendar}
\address{Isabelle Chalendar,  Universit\'e Gustave Eiffel, LAMA, (UMR 8050), UPEM, UPEC, CNRS, F-77454, Marne-la-Vallée (France)}
\email{isabelle.chalendar@univ-eiffel.fr}

\author{J.R. Partington}
\address{Jonathan R. Partington, School of Mathematics, University of Leeds, Leeds LS2 9JT, UK}
 \email{J.R.Partington@leeds.ac.uk}
%
%
%
%
%
%
\begin{abstract}
This paper considers composition operators on Zen spaces (a class of weighted Bergman spaces of the
right half-plane related to weighted function spaces on the positive half-line by means of the
Laplace transform). Generalizations are given to work of Kucik on norms and essential norms, to work
of Schroderus on (essential) spectra, and to work by Arvanitidis and the authors on semigroups of
composition operators. 
The results are illustrated by consideration of the Hardy--Bergman space; that is, the intersection of the
Hardy  and Bergman Hilbert spaces on the half-plane.
\end{abstract}

 \maketitle
 
 {\bf MSC class:} {30H10, 30H20, 47B33, 47D03}

{\bf Keywords:} {Composition operator, Hardy space, Bergman space, spectrum, essential spectrum, operator semigroup}

\section{Introduction and background material}

The subject of this paper involves properties of composition operators 
on holomorphic function spaces on the right half-plane $\CC_+$, both
as individual operators and as elements of one-parameter semigroups.
One difficulty, even in the case of the Hardy space $H^2(\CC_+)$, is that not all
composition operators on the spaces are bounded, although in many cases the
bounded operators have been characterised, as we explain below.
The literature on   semigroups of composition operators is not as extensive as it is for the 
disc, although it has been an object of study since the work of Berkson and Porta \cite{BP78}.

In this note we shall concentrate on the so-called Zen spaces (weighted Hardy--Bergman spaces),
which form a large class of spaces with applications in systems and control theory
\cite{JPP13,JPP14}. For these
information on the bounded composition operators is also available, thanks to Kucik \cite{kucik18}.

The remainder of this section presents necessary background material.
In Section \ref{sec:2} we provide new results on the norm and spectral radius of composition operators.
Then, in Section \ref{sec:spectral}, we
consider the spectral theory of composition operators on the half-plane with linear fractional symbols, extending several
results of Schroderus \cite{schroderus}.
Finally, in Section \ref{sec:4}, we consider semigroups of composition operators on Zen spaces, providing generalizations
of results of Arvanitidis \cite{arv15}.

\subsection{Zen spaces}

Kucik \cite{kucik18} considered composition operators on Zen spaces $A^2_\nu$,
which are isometrically Laplace transforms of weighted Hardy spaces $L^2_w(0,\infty)=L^2(0,\infty, w(t) dt)$.
For general background see   \cite{CP14}, but we provide the basic facts now.

Let $\nu$ be a positive regular Borel measure on $[0,\infty)$ satisfying the doubling condition
\[
R:= \sup_{t>0} \frac{\nu[0,2t)}{\nu[0,t)} < \infty.
\]
The Zen space $A^2_\nu$ is defined to consist of all analytic  functions $F$ on $\CC_+$
such that the norm, given by
\[
\|F\|^2=\sup_{\epsilon>0} \int_{\overline{\CC_+}} |F(s+\epsilon)|^2 \, d\nu(x) \, dy
\]
is finite, where we write $s=x+iy$ for $x\ge 0$ and $y \in \RR$.

The best-known examples here are:
\begin{enumerate}
\item For  $\nu=\delta_0$, a Dirac mass at $0$,  we obtain the
Hardy space $H^2(\CC_+)$;
\item For $\nu$ equal to Lebesgue measure  we obtain the Bergman space $A^2(\CC_+)$.
\end{enumerate}

Often we shall have $\nu(\{0\})=0$, in which case $\|F\|^2$ can be written simply as
\[
\int_{\overline{\CC_+}} |F(s)|^2 \, d\nu(x) \, dy.
\]

\begin{thm}\cite{JPP13}
Suppose that $w$ is given as a weighted Laplace transform
\[
w(t)=2\pi \int_0^\infty e^{-2rt} \, d\nu(r), \qquad (t >0).
\]
Then the Laplace transform provides an isometric map 
\[
\LL:  L^2(0,\infty,   w(t) dt ) \to A^2_\nu.  
\]
\end{thm} 

The result also holds for Hilbert-space valued functions \cite{AP21}.\\

The following result   generalizes the Elliott--Jury theorem in \cite{EJ12} and the Elliott--Wynn theorem
in \cite{EW11}.

\begin{thm}\label{thm:normcphi}\cite{kucik18}
The composition operator $C_\phi$ is bounded on $A^2_\nu$ if and only if
$\phi$ has a finite nonzero angular derivative $L=\angle\lim_{z \to \infty} z/\phi(z)$ at $\infty$.
In that case
\beq\label{eq:opnorm}
L \inf_{t>0} \frac{w(t)}{w(L t)} \le \|C_\phi\|^2 \le  L \sup_{t>0} \frac{w(t)}{w( L t)}.
\eeq
\end{thm}
This gives correctly $\|C_\phi\|=\sqrt L$ in the Hardy case and $\|C_\phi\|=L$ in the 
Bergman case.

\subsection{Semigroups of composition operators}

Berkson and Porta \cite{BP78} give the following criterion for
an analytic function $G$ to generate a one-parameter semigroup
of analytic mappings of the right half-plane $\CC_+$ to itself; that is, solutions to
\beq\label{eq:phig}
\frac{\partial \phi_t(z)}{\partial t}=G(\phi_t(z)), \qquad \phi_0(z)=z,
\eeq
in terms of the following condition:
\beq\label{eq:nscbp}
x \frac{\partial(\re G)}{\partial x} \le \re G \qquad \hbox{on} \quad \CC_+,
\eeq
where $x=\re z$. The associated  composition operators $C_{\phi_t}$ are bounded on $H^2(\CC_+)$  if and only if
the non-tangential limit
\[
L_t:=\angle \lim_{z \to \infty}z/\phi_t(z)
\]
 exists and is non-zero. It will be positive, and then $\|C_{\phi_t}\|=L_t^{1/2}$ (see \cite{EJ12,matache}).

Arvanitidis \cite{arv15} showed that a necessary and sufficient condition for boundedness of the composition operators
is that $\delta:=\angle\lim_{z \to \infty} G(z)/z$ exists. In this case $\|C_{\phi_t}\|=e^{-\delta t/2}$ and
the semigroup is quasicontractive.

The following result is taken from   \cite{ACP16}.
\begin{thm}
For an operator $A$ given by $Af=Gf'$ on $D(A) \subseteq H^2(\CC_+)$, the following are equivalent:\\
(i) $A$ generates a quasi-contractive $C_0$-semigroup of bounded composition operators on $H^2(\CC_+)$;\\
(ii) Condition \eqref{eq:nscbp} holds and  $\angle \lim_{z \to \infty} G(z)/z$ exists.
\end{thm}

\section{The essential norm and spectral radius}\label{sec:2}

We know from \cite[Lem.~4]{kucik18} that in every $A^2_\nu$ space the normalized reproducing kernels $k_n / \|k_n\|$ tend weakly to $0$
as $n \to\infty$. Now we can adapt the proof that there is no compact composition
operator as follows.

Using the isometry between $L^2(0,\infty; w(t) \, dt)$ and $A^2_\nu$ we may pull back the reproducing kernel
of $A^2_\nu$ at a point $\lambda \in \CC_+$ 
by the inverse Laplace transform to obtain
$k_\lambda(t):=e^{-\overline\lambda t}/w(t)$ and its norm is $\langle k_\lambda,k_\lambda\rangle^{1/2}$ or
\beq\label{eq:wint}
\left(\int_0^\infty \frac{e^{-2\re \lambda t}}{w(t)} \, dt \right)^{1/2}.
\eeq

\begin{thm}
Let $C_\phi$ be a bounded composition operator on the Zen space $A^2_\nu$ corresponding to
a weight $w$ on $(0,\infty)$. Let $L$ denote the 
 finite nonzero angular derivative $L=\angle\lim_{z \to \infty} z/\phi(z)$ at $\infty$.
Then
\beq\label{eq:essnorm}
L \inf_{t>0} \frac{w(t)}{w(L t)} \le \|C_\phi\|_e^2 \le \|C_\phi\|^2 \le  L \sup_{t>0} \frac{w(t)}{w( L t)}.
\eeq
\end{thm}
\beginpf
For every $\delta>0$ there is a compact operator $Q$ on $A^2_\nu$ such that 
\begin{align*}
\|C_\phi\|_e + \delta &\ge 
\|C_\phi-Q\| \\ &\ge \limsup_{n \to \infty} \|(C_\phi-Q)^*k_n\|/\|k_n\| \\ &=  \limsup_{n \to \infty} \|k_{\phi(n)}\|/\|k_n\|
\end{align*}
since $Q^* k_n / \|k_n\| \to 0$ in norm. Since this is true for all $\delta>0$ we have
\[
\|C_\phi\|^2_e \ge \limsup_{n \to \infty} \left(\int_0^\infty \frac{e^{-2\re \phi(n) t}}{w(t)} \, dt \right) /
\left(\int_0^\infty \frac{e^{-2  n t}}{w(t)} \, dt \right)
\]
by \eqref{eq:wint}. Hence, for every $0<M<L$ we have
\begin{align*}
\|C_\phi\|^2_e &\ge \limsup_{n \to \infty} \left(\int_0^\infty \frac{e^{-2 n t/M}}{w(t)} \, dt \right) /
\left(\int_0^\infty \frac{e^{-2  n t}}{w(t)} \, dt \right)\\
&= \limsup_{n \to \infty} \left(\int_0^\infty \frac{e^{-2 un}}{w( Mu)}  M \,du \right) /
\left(\int_0^\infty \frac{e^{-2  n t}}{w(t)} \, dt \right)\\
& \ge M \inf_{u > 0} \frac{ w(u)}{ w(Mu)}.
\end{align*}

Now, using Theorem \ref{thm:normcphi}, we have the estimate for the essential norm given in
\eqref{eq:essnorm}.
\endpf

For the case of the Hardy and weighted Bergman spaces we therefore recover the formula $\|C_\phi\|=\|C_\phi\|_e$ 
from \cite{EJ12,EW11}.\\

Note that the left and right hand sides of \eqref{eq:opnorm} and \eqref{eq:essnorm} may certainly differ, and 
as an example we shall consider the
Hardy--Bergman space $H^2(\CC_+) \cap A^2(\CC_+)$.
On the disc the Hardy space $H^2(\DD)$ is contained in the Bergman space $A^2(\DD)$
but no such inclusion either way exists on the half-plane, as we shall now explain.
We may give  $H^2(\CC_+) \cap A^2(\CC_+)$ a Hilbert space structure with the norm defined by
\[
\|f\|^2_{H^2(\CC_+) \cap A^2(\CC_+)} = \|f\|^2_{H^2(\CC_+)}+\|f\|^2_{A^2(\CC_+)}.
\]
Since  the Hardy space is isomorphic by means of the Laplace transform to
$L^2(0,\infty)$ and the Bergman space to $L^2(0, \infty, \, dt/t)$ we see
that neither space is contained in the other and that the Hardy--Bergman
space $H^2(\CC_+) \cap A^2(\CC_+)$ corresponds to $L^2(0,\infty,w(t)\, dt )$
with the
weight $w(t)=1+1/t$.

In this  case $w(t)/w(Lt)=\dfrac{t+1}{t+L}$, so that the sup and inf are
$1$ and $1/L$ in some order.\\

For the spectral radius $\rho(C_\phi)$ we may apply \eqref{eq:opnorm} to get an estimate which is
once more tight in the case of the Hardy and weighted Bergman spaces.

\begin{cor}
For a bounded composition $C_\phi$ on a space $A^2_\nu$ corresponding to a weight $w$, the spectral radius $\rho(C_\phi)$
satisfies
\[
L\limsup_{n \to \infty}\left( \inf_{t>0} \frac{w(t)}{w(L^n t)}\right)^{1/n} \le \|\rho(C_\phi)\|^2 \le  L 
\liminf_{n \to \infty}\left(\sup_{t>0} \frac{w(t)}{w( L^n t)}\right)^{1/n},
\]
where $L=\angle\lim_{z \to \infty} z/\phi(z)$.
\end{cor}
\beginpf
For the iterate $\phi^{(n)}$ the corresponding 
angular derivative $L_n=\angle\lim_{z \to \infty} z/\phi^{(n)}(z)$ is simply $L^n$, being the product of $n$ terms 
each tending to $L$. Thus 
\[
L^n \inf_{t>0} \frac{w(t)}{w(L^n t)} \le \|C_\phi^n\|^2 \le  L^n \sup_{t>0} \frac{w(t)}{w( L^n t)},
\]
from which the result follows by the standard spectral radius formula.
\endpf

\section{Spectral theory}\label{sec:spectral}

For the spectral theory of composition operators on the half-plane, little
seems to be known, although Schroderus \cite{schroderus} has determined
the spectrum and essential spectrum in the case of linear fractional mappings $\phi$ for Hardy and weighted Bergman
spaces. The associated composition operator $C_\phi$ is bounded if and only if $\phi$
is a parabolic or hyperbolic mapping fixing $\infty$; that is,
$\phi(s)=\mu s+s_0$ with $\mu>0$ and $\re s_0 \ge 0$ (we shall discuss this in detail below).
The same applies in general Zen spaces, by Theorem \ref{thm:normcphi}.

Schroderus restricts herself to Hardy spaces and weighted Bergman spaces on the upper half-plane.
Transforming to the right half-plane these are Zen spaces $A^2_\nu$ with measures $d\nu(x)=x^\alpha \,dx$ for $\alpha>-1$ (in the following the Hardy space is formally identified with the case $\alpha=-1$). Her results in
this particular context are:

\begin{thm}\cite[Thm. 1.1, Thm 1.2]{schroderus}\label{thm:riikka}
\begin{enumerate}
\item In the parabolic case $\phi(s)=s+s_0$, where $\re s_0 \ge 0$ and $s_0 \ne 0$, the spectrum and essential spectrum
of $C_\phi$ coincide and equal\\
(a) $\TT$ if $s_0 \in i\RR$;\\
(b) $\{e^{-s_0 t}: t \ge 0 \} \cup \{0\}$ if $s_0 \in \CC_+$.
\item In the hyperbolic case $\phi(s)=\mu s + s_0$, where $\mu \in (0,1) \cup (1,\infty)$ and
$\re s_0 \ge 0$, the spectrum and essential spectrum
of $C_\phi$ coincide and equal\\
(a) $\{ \lambda \in \CC: |\lambda| = \mu^{-(\alpha+2)/2} \}$ if $s_0 \in i\RR$;\\
(b)  $\{ \lambda \in \CC: |\lambda| \le \mu^{-(\alpha+2)/2} \}$ if $s_0 \in \CC_+$.\\
\end{enumerate}
\end{thm}

\subsection{The parabolic case}

For a general Zen space, part (1) of Theorem \ref{thm:riikka} is easy to prove.

\begin{prop}
Let $\phi(s)=s+s_0$ where $\re s_0 \ge 0$ and $s_0 \ne 0$. Then the
spectrum and essential spectrum of $C_\phi$ on $A^2_\nu$
coincide and equal\\
(a) $\TT$ if $s_0 \in i\RR$;\\
(b) $\{e^{-s_0 t}: t \ge 0 \} \cup \{0\}$ if $s_0 \in \CC_+$.
\end{prop}
\beginpf
The composition operator is seen to be unitarily equivalent to
the multiplication operator on $L^2(0,\infty,   w(t) dt )$
given by multiplication by the function $t \mapsto e^{-s_0 t}$.
Thus it is a normal operator: its spectrum and essential spectrum
equal the closure of $\{ e^{-s_0 t}: t \ge 0 \}$,  as in (a) and (b).
\endpf

\subsection{The hyperbolic case}

In this subsection we take
$\phi(s)=\mu s + s_0$, where $\mu \in (0,1) \cup (1,\infty)$ and
$\re s_0 \ge 0$.

Schroderus observed that we have the following simplifications:

For $s_0=iy$, with $y \in \RR$, the operator $C_\phi$ is similar to $C_\psi$ with $\psi(s)=\mu s$.
Indeed, if $\rho(s)=s+iy/(\mu-1)$ then
\[
\rho^{-1} \circ \psi \circ \rho (s) = \mu(s+iy/(\mu-1))-iy/(\mu-1) = \mu s + iy.
\]

Similarly, for $s_0=x+iy$ with $x>0$ and $y \in \RR$, let
$\psi(s)=\mu s +   x$ and take the same $\rho$. Then
\[
\rho^{-1} \circ \psi \circ \rho (s) = \mu(s+iy/(\mu-1))+x-iy/(\mu-1) = \mu s + x+iy.
\]
Thus, since $C_\rho$ is a unitary map on every $A^2_\nu$, we need only consider the spectrum of $C_\phi$
when $\phi(s)=\mu s+ x$ for $\mu \in (0,1) \cup (1,\infty)$ and $x \ge 0$.\\

In case (a) we can estimate the spectral radius of $C_\phi$ and $C^{-1}_\phi$,
although this tells us only that the spectrum is contained in an annulus (which may be a circle).
In case (b) we have that the spectrum is contained in a disc.\\

We showed in \cite[Prop. 3.5]{CP14} that the norm of the weighted composition operator $C_\psi$ 
corresponding to $\psi(s)=\mu s$ with $\mu >0$ is exactly $\ds \left(\frac{1}{\mu} \sup_{t >0} \frac {w(\mu t)}{w(t)} \right)^{1/2}$.
The method of \cite{CP14} also shows that for  all $x \ge 0$
and $\phi(s)=\mu s+ x$ the composition operator has norm equal to the norm of
the mapping $f \mapsto g$ in $L^2(0,\infty; w(t) \, dt)$, where
\[
g(t) = \frac{1}{\mu} e^{-x t/\mu}   f(t/\mu)
\]
since
\begin{align*}
\int_0^\infty e^{-st} g(t) \, dt &=\frac{1}{\mu}\int_0^\infty e^{-st}   e^{-x t/\mu}   f(t/\mu) \, dt\\
&= \int_0^\infty e^{-\mu s\tau } e^{-x\tau}   f(\tau) \, d\tau
\end{align*}
and 
\[
\|g\|^2 = \frac{1}{\mu^2} \int_0^\infty e^{-2x t/\mu} |f(t/\mu)|^2 w(t) \, dt = \frac{1}{\mu} \int_0^\infty e^{-2x \tau}
|f(\tau)|^2 w(\mu \tau) \, d\tau,
\]
whence
\[
\|C_\phi\|^2 = \frac{1}{\mu} \sup_{t>0} \frac{e^{-2x t}
 w(\mu t)}
{w(t)}.
\]
Now by looking at iterates we can obtain an explicit formula for the spectral radius of $C_\phi$ as
the $n$th iterate of $\phi$ is given by
$
\phi_n(s)= \mu^n s + x_n$, where $x_n=\dfrac{\mu^n-1}{\mu-1}x$, and this has the same form as $\phi$.\\

That is,
\beq\label{eq:specrad}
\rho(C_\phi)= \frac{1}{\sqrt \mu}\lim_{n \to \infty}  \left( \sup_{t>0} \frac{e^{-2x_n t}
 w(\mu^n t)}
{w(t)} \right)^{1/2n}.
\eeq

In the time domain (that is, using the inverse Laplace transform), the operator $C_\phi$ corresponds to the
operator  $B_\phi$ on $L^2(0,\infty, w(t) dt)$ with
\[ B_\phi f(t) = af(bt)e^{-ct} \]
such that we have
\begin{align*}
\int_0^\infty a f(bt) e^{-ct} e^{-st} dt &= \int_0^\infty a f(u) e^{-cu/b} e^{-su/b} du/b \\
& = (a/b) (\LL f)(s/b + c/b) = (\LL f)(\mu s+ x),
\end{align*}
and so $a=b=1/\mu $ and $c=x/\mu$.\\

With $Af(t)= \frac{1}{\mu} f(t/\mu)e^{-xt/\mu}$ on $L^2(0,\infty,   w(t) dt )$ we find that
\[
\int_0^\infty \frac{1}{\mu} f(t/\mu) e^{-xt/\mu}\overline{g(t)}w(t)\, dt
=
\int_0^\infty f(u)  e^{-xu}\overline{g(\mu u)} \frac{w(\mu u)}{w(u)} w(u) \, du
\]
so that
$\ds A^*g(u)= g(\mu u) e^{-xu}\frac{w(\mu u)}{w(u)}$.

\subsubsection{The case $x=0$}

It is enough to consider the spectrum of $C_\phi$ for $\phi(s)=\mu s$ with $\mu > 1$
since the case $0<\mu < 1$ may be studied by taking inverses.

\begin{thm}
For the composition operator $C_\phi$ on $A^2_\nu$ where
$\phi(s)=\mu s$ with $\mu>1$ and $\nu$ determining the weight $w$ on $(0,\infty)$
we have
\[
\sigma(C_\phi) \subseteq \{z \in \CC: r \le |z| \le R \},
\]
where 
\[
r=\frac{1}{\sqrt\mu} \lim_{n \to \infty}\left(  \inf_{t >0}  \frac{ w(\mu^n t)}
{w(t)}\right)^{1/(2n)}
\]
 and 
\[
R=\frac{1}{\sqrt\mu} \lim_{n \to \infty}\left(  \sup_{t >0}  \frac{ w(\mu^n t)}
{w(t)}\right)^{1/(2n)}.
\]
\end{thm}
\beginpf
This follows from
\eqref{eq:specrad}.
\endpf

With $\ds A^*f(t)= f(\mu t)  \frac{w(\mu t)}{w(t)}$,
we have
 eigenvalues of $A^*$ given by
 $\mu^\alpha$, with eigenvectors $f(t)=t^\alpha/w(t)$,
 provided that these eigenvectors lie in the space $L^2(0,\infty,w(t) dt)$.
 Indeed the eigenvalues then have infinite multiplicity, since
 $f\chi_E$ is an eigenvector for any measurable $E \subset \RR_+$ such that $\mu E=E$,   the notation $\chi_E$ denoting  
the characteristic (indicator) function of a set $E$.
Clearly there are infinitely many distinct subsets $E \subset \RR$ such that $\mu E=E$. For example,
if $\mu>1$ we may take 
\[
E=\bigcup_{n=-\infty}^\infty \chi_{(\mu^n,\mu^n (1+\delta))}
\]
for any $\delta$ with $0 < \delta < \mu-1$. 
The same is true for $0<\mu<1$ working with $1/\mu$ instead of $\mu$: this will be required later.
 
 Thus we have the following lower bound for the essential spectrum $\sigma_e(C_\phi)$:
 \begin{thm}\label{thm:eigenspace}
For the composition operator $C_\phi$ on $A^2_\nu$ where
$\phi(s)=\mu s$ with $\mu>1$ and $\nu$ determining the weight $w$ on $(0,\infty)$
it holds that $\sigma_e(C_\phi)$ contains
all $\alpha$ such that
\beq\label{eq:eigenspace}
 \int_0^\infty |t^\alpha|^2/w(t) < \infty.
 \eeq
\end{thm}
\beginpf
The eigenvectors of $f(t)=t^\alpha/w(t)$ of the unitarily equivalent operator $A^*$ lie in the space  $L^2(0,\infty,w(t) dt)$
if and only if 
\eqref{eq:eigenspace} holds.
\endpf

Note that the $\alpha$ occurring in Theorem \ref{thm:eigenspace}
form an annulus (if the set is non-empty), and in some cases (e.g. the Hardy--Bergman space below)
this enables us to find the spectrum exactly.

\subsubsection{The case $x>0$} 

Here the calculations are necessarily more difficult, but we do have the
spectral radius formula \eqref{eq:specrad} to give a disc in which the spectrum is contained.
In many cases all the points in the interior are eigenvalues, either of $A$ or $A^*$.
Indeed, they are again eigenvalues of infinite multiplicity, and hence in the essential spectrum,
as we see in the following two propositions,

\begin{prop}\label{prop:spec1}
For $\mu>1$, if $\alpha \in \CC$ and $\beta=x/(\mu-1)>0$ are such that
the function $f:t \mapsto t^\alpha e^{-\beta t}$ lies in $L^2(0,\infty, w(t) dt)$,   then
$A$ has the eigenvalue $1/\mu^{\alpha+1}$ 
and $f\chi_E $  is an eigenvector for any measurable $E \subset \RR_+$ such that $\mu E=E$.
\end{prop}
\beginpf
This is an easy calculation since
$Af(t)= \frac{1}{\mu} f(t/\mu)e^{-xt/\mu}$.
\endpf

\begin{prop}\label{prop:spec2}
For $0<\mu<1$ if $\alpha \in \CC$ and $\beta=x/(1-\mu)>0$ are such that
the function $f(t)=t^\alpha e^{-\beta t}/w(t)$ lies in $L^2(0,\infty, w(t) dt)$,   then
$A^*$ has the eigenvalue $\mu^\alpha$ 
and $f\chi_E$ is an eigenvector for any measurable $E \subset \RR_+$ such that $\mu E=E$.
\end{prop}
\beginpf
This is similarly straightforward, since
\[
A^*f(t)= f(\mu t) e^{-xt}  {w(\mu t)}/{w(t)}.
\]
\endpf

One technique from \cite{schroderus} is not available in general: for the Hardy and standard
weighted Bergman spaces, the weighted composition operators with $\mu$ and $1/\mu$
are related by taking adjoints. This is not true unless  $w(t)$ is a power of $t$.

However, as we shall see in the next section the information here is enough to allow us to
determine the exact spectrum in an important new example (as well as the cases covered by Schroderus).

\subsection{The Hardy--Bergman space}

As a significant example of a situation not covered by older results, 
let us look again at the Hardy--Bergman space $H^2(\CC_+) \cap A^2(\CC_+)$. To
determine the spectrum we may take any equivalent norm, so we take $w(t)=1+1/t$ as before.
\begin{thm}\label{th:hh}
	Let $X=H^2(\CC_+) \cap A^2(\CC_+)$ be the Hardy-Bergman space and $C_\varphi\in {\mathcal L}(X)$. 
	\begin{enumerate}
		\item[a)] If $\varphi(s)=\mu s$ with $\mu>0$, then $\sigma(C_\varphi)$ is the annulus 
		\[
		\sigma(C_\phi)= \left\{z \in \CC:\min\left \{ \frac{1}{\mu}, \frac{1}{\sqrt \mu}\right\} \le |z| \le 
		 \max\left \{ \frac{1}{\mu}, \frac{1}{\sqrt \mu}\right\}\right\}.
		\]
		\item[b)] If $\varphi(s)=\mu s + (x+iy)$ with $\mu,x>0$ and $y\in\RR$, then $\sigma(C_\varphi)$ is the closed disc 
		\[   \sigma(C_\phi)= \{ z \in \CC:   |z| \le 1/\mu\}.\] 
	\end{enumerate} 
	In all cases we have $\sigma_e(C_\phi)=\sigma(C_\phi)$.
\end{thm}
The proof of this theorem follows from the two subsections detailed below. 
\subsubsection{The case $x=0$}
For $\phi(s)=\mu s$ with $\mu>0$ we have
\[
\|C_\phi\|^2 =  \frac{1}{\mu} \sup_{t>0}  \frac{ 
 w(\mu t)}
{w(t)} =  \frac{1}{\mu}  \sup_{t>0} \frac{1+1/(\mu t)}{1+1/t}= \frac{1}{\mu}  \sup_{t>0} \frac{\mu t+1 }{\mu t+\mu}.
\]
This is 
\[
\begin{cases}
1/\mu^2 & \hbox{if} \quad 0<\mu<1,\\
1/\mu & \hbox{if} \quad  \mu>1. 
\end{cases}
\]
Thus for $\mu>1$ we have
\[
\|C^n_\phi\|^2 =  \frac{1}{\mu^n}
\]
while
\[
\|C^{-n}_\phi\|^2  = \mu^{2n}.
\]
That is, the spectrum of $C_\phi$ lies in the annulus $\{z \in \CC: \frac{1}{\mu} \le |z| \le \frac{1}{\sqrt \mu}\}$,
but not in a circle. So we see that Theorem \ref{thm:riikka} (2)(a) does not hold in the general case of a Zen space.

We now show that the (essential) spectrum equals the whole annulus. To do this we use Theorem \ref{thm:eigenspace}.
The function $t^\alpha/(1+1/t)$ lies in the space $L^2(0,\infty, (1+1/t) \, dt)$ if and only if
\[
\int_0^\infty \frac{|t^\alpha|^2}{(1+1/t)^2} (1+1/t) \, dt < \infty,
\]
 so that $-1<\re \alpha < -1/2$. From this we see that all points in the annulus
 $\{z \in \CC: \frac{1}{\mu} < |z| < \frac{1}{\sqrt \mu}\}$ are eigenvalues of $A^*$ so 
 that for $\mu>1$ we have
 \[
 \sigma(C_\phi)= \left\{z \in \CC: \frac{1}{\mu} \le |z| \le \frac{1}{\sqrt \mu}\right\}.
 \]
 
 By considering inverses, we see that for $0<\mu<1$
  \[
 \sigma(C_\phi)= \left\{z \in \CC: \frac{1}{\sqrt \mu} \le |z| \le \frac{1}{\mu}\right\}.
 \]

 \subsubsection{The case $x>0$}
We begin with the spectral radius formula \eqref{eq:specrad}.
Elementary calculus shows that the supremum in \eqref{eq:specrad} is at $t=0$ except if
$2x_n < \mu^n-1$. This requires $\mu>1$.
If the supremum is  at $t=0$ it is $1/\mu^n$.\\

Otherwise (with $\mu>1$) we can estimate the supremum for each $n$ by considering the two cases
(recall that $x_n \to \infty$ as $n \to \infty$):
\begin{enumerate}
\item  $  t \ge  \dfrac{n\log \mu}{2x_n }$,
when the expression is clearly at most $\mu^{-n}$ (the second factor is bounded by 1);\\
\item $0 \le   t \le   \dfrac{n\log \mu}{2x_n }$, when we may ignore the exponential and obtain an upper bound 
\[
\frac{1+1/(\mu^n t)}{1+1/ t} \le \frac{1+2x_n/(n\mu^n \log \mu)}{2x_n/n \log \mu} \le \frac{1+C_1/n}{C_2 \mu^{n-1}/n},
\]
where $C_1$ and $C_2$ are independent of $n$.
\end{enumerate}

Taking the maximum of these two estimates and then the $2n$-th root gives a limit $1/\sqrt{\mu}$.
That is $\rho(C_\phi) \le \frac{1}{\sqrt \mu} \times \frac{1}{\sqrt \mu} = \frac{1}{\mu}$
for all $x>0$.\\

Now we use Propositions \ref{prop:spec1} and \ref{prop:spec2} to find eigenvectors. 

For $\mu>1$
the range of admissible $\alpha$ is $\{\re \alpha>0\}$ and the eigenvalue is $1/\mu^{\alpha+1}$, so
that every point of the disc $\{ z \in \CC: 0 < |z| < 1/\mu\}$ is an eigenvector of $A$ and
we deduce that
$\sigma(C_\phi)= \{ z \in \CC:   |z| \le 1/\mu\}$.

Similarly, for $0< \mu < 1$ the range of admissible $\alpha$ is $\{\re \alpha>-1\}$ and the eigenvalue
is $\mu^\alpha$, so again we deduce that 
$\sigma(C_\phi)= \{ z \in \CC:   |z| \le 1/\mu\}$, and the same holds for $\sigma_e(C_\phi)$.

\section{Semigroups of composition operators}
\label{sec:4}

\subsection{Norm estimates}

Our first task is to generalize
the result of
Arvanitidis \cite{arv15} showing that a necessary and sufficient condition for boundedness of
a semigroup $(C_{\phi_t})$ of composition operators
on $H^2(\CC_+)$ with 
generator $A:f \mapsto Gf'$
is that $\delta:=\angle\lim_{z \to \infty} G(z)/z$ exists. In this case $\|C_{\phi_t}\|_{H^2(\CC_+)}=e^{-\delta t/2}$.
This relies on a theorem in \cite{CDMP06} (which is stated for the disc), which gives the angular derivative,
and hence the norm, in terms of $G$. In particular, $L_t=e^{-\delta t}$. 

By virtue of Theorem \ref{thm:normcphi} (Kucik) we know that the condition for boundedness does not depend on
which Zen space we use, and all that changes is the norm of the operator. We may thus
state the following easy theorem.

\begin{thm}
For a semigroup of analytic self-maps $(\phi_t)$ on $\CC_+$, the following are equivalent:
\begin{enumerate}
\item the non-tangential limit $\delta:=\angle\lim_{z \to \infty} G(z)/z$ exists;
\item the semigroup $(C_{\phi_t})_{t \ge 0}$ consists of bounded operators on $A^2_\nu$.
\end{enumerate}
In this case, we have 
\[
L _t\inf_{x>0} \frac{w(x)}{w(L_t x)} \le \|C_{\phi_t}\|^2 \le  L_t \sup_{x>0} \frac{w(x)}{w( L_t x)},
\]
where $L_t=e^{-\delta t}$ for $t \ge 0$.
\end{thm}
\beginpf
This follows from 
\cite[Thm. 3.4]{arv15} together with the estimate 
 \eqref{eq:opnorm} from \cite{kucik18}.
 \endpf
 
 For example, with the 
 Hardy space $H^2(\CC_+)$ considered by Arvanitidis, we have
 $\|C_{\phi_t}\|=e^{-\delta t/2}$, while for the 
 standard Bergman space $A^2(\CC)$ we have the apparently new result that 
$\|C_{\phi_t}\|=e^{-\delta t}$.

In \cite{AC19} there is a result applying to analytic function spaces, i.e., Banach spaces $X$ of holomorphic functions on a domain $\Omega$ for which
point evaluations are continuous functionals. The result applies to such spaces satisfying the following supplementary
condition:

(E) \quad 
If $(z_n)$ is a sequence in $\Omega$ such that $z_n \to z \in \overline\Omega \cup \{\infty\}$ and $\lim_{n \to \infty} f(z_n)$ exists in $\CC$  for all $f \in X$
then $z \in \Omega$.

\begin{thm}\label{thm:E}
Suppose that $(T_t)_{t \ge 0}$ is a $C_0$ semigroup on a function space $X$ with property (E) such that 
 for some $G \in \Hol(\Omega)$ the generator $A$
is the operator $f \mapsto Gf'$ for all $f \in \dom(A)$. Then $(T_t)$ is a semigroup of composition operators.
\end{thm}

Similar results are given in \cite{GY19}, for function spaces on the disc $\DD$.
We may now apply Theorem \ref{thm:E} to our situation.

\begin{thm}
Suppose that $(T_t)_{t \ge 0}$ is a $C_0$ semigroup on a Zen space $A^2_\nu$
such that for some $G \in \Hol(\CC_+)$ the generator $A$
is the operator $f \mapsto Gf'$ for all $f \in \dom(A)$. Then $(T_t)$ is a semigroup of composition operators.
\end{thm}

\beginpf
For the halfplane we need to find sequences $(z_n)$ in $\CC_+$ tending to $\infty$ or a point on the imaginary axis for which 
$\lim f(z_n)$ does not exist for some $f \in A^2_\nu$.
 

Now if $(z_n)$ is such that $(f(z_n))$ converges for all $f \in A^2_\nu$ then by the uniform boundedness theorem
the reproducing kernels at $(z_n)$ are uniformly bounded. But if $\re z_n \to 0$ monotonically then the
integrals in \eqref{eq:wint} with $\lambda=z_n$ tend to $\infty$ (note that $w$ is a decreasing function).
This shows that condition (E) holds for $A^2_\nu$ (take sequences tending to either a point on the imaginary axis or to infinity with real parts tending to 0). The result follows from Theorem \ref{thm:E}.
\endpf

\begin{rem}{\rm
We may also consider the question of groups of composition operators on $A^2_\nu$.
However, since we have same semigroups for all the spaces under consideration, with different norms, we have the same groups. 
Proposition 4.4 of \cite{ACP16}
asserts that for a $C_0$-quasicontractive group of bounded composition operators on on $H^2(\CC_+)$ we have
$G(z)=ps+iq$ for real $p$ and $q$, and $\phi_t(s)=e^{pt}s+\frac{iq}{p}(e^{pt}-1)$
if $p \ne 0$, or $\phi_t(s)=s+iqt$ if $p=0$.
The same will apply to any Zen space $A^2_\nu$.
}
\end{rem}

\subsection{Semigroups of linear fractional mappings}

As in Section \ref{sec:spectral} we may discuss parabolic or hyperbolic mappings fixing $\infty$.
That is, $\phi_t(s)=\mu_t s + s_t$, where $\mu_t > 0$ and $s_t \in \overline{\CC_+}$.

The semigroup relation $\phi_{t+u}=\phi_t \circ \phi_u$ gives us $\mu_t=e^{pt}$ for some $p \in \RR$, and
\begin{enumerate}
\item $s_t= \alpha (e^{pt}-1)$ for some $\alpha \in \overline{\CC_+}$ if $p \ne 0$ (the hyperbolic case), or
\item $s_t=\alpha t$ for some $\alpha \in \overline{\CC_+}$ if $p=0$ (the parabolic case).
\end{enumerate}
Using \eqref{eq:phig} we see that $G(z)=pz+p\alpha$ in the hyperbolic case and $G(z)=\alpha$ in the parabolic case.

An immediate corollary of Theorem~\ref{th:hh} is the following description of the spectrum.   
\begin{cor}\label{cor:hh}
	Let $X=H^2(\CC_+) \cap A^2(\CC_+)$ be the Hardy-Bergman space and $(C_{\varphi_t})_{t\geq 0}\subset {\mathcal L}(X)$. . 
	\begin{enumerate}
		\item[a)] If $\varphi_t(s)=e^{pt} s + \alpha(e^{pt}-1)$ with $p\neq 0$, then $\sigma(C_{\varphi_t})$ is the annulus 
		\[
		\sigma(C_{\phi_t})= \left\{z \in \CC:\min\left \{ \frac{1}{e^{pt}}, \frac{1}{e^{pt/2}}\right\} \le |z| \le 
		\max\left \{ \frac{1}{e^{pt}}, \frac{1}{e^{pt/2}}\right\}\right\}.
		\]
		\item[b)] If $\varphi_t(s)=s + \alpha t$, then $\sigma(C_{\varphi_t})$ is the closed disc 
		\[   \sigma(C_{\phi_t})= \{ z \in \CC:   |z| \le e^{-pt}\}\] 
	\end{enumerate} 
	
	In all cases the spectrum and essential spectrum coincide.
\end{cor}

\end{document}